\documentclass{amsart}
\usepackage[english]{babel}
\usepackage{amsfonts,amssymb,amsmath}

\newtheorem{thm}{Theorem}
\newtheorem*{quest}{Question}

\theoremstyle{remark}
\newtheorem*{rem}{Remark}
\newtheorem*{ack}{Acknowledgements}

\newcommand{\de}{\partial}

\newcommand{\db}{\overline{\partial}}
\newcommand{\Ric}{\mathrm{Ric}}

\begin{document}
\title[K\"AHLER GEOMETRY]{On the critical points of the $E_k$ functionals in K\"ahler geometry}
\author{{Valentino Tosatti}}%
\address{Department of Mathematics, Harvard University, 1 Oxford St, Cambridge MA 02138}%
\email{tosatti@math.harvard.edu}%
\thanks{The author is supported by a Harvard Mathematics Department grant.}
\begin{abstract}
We prove that a K\"ahler metric in the anticanonical class which is a critical point of the functional $E_k$ and
has nonnegative Ricci curvature, is necessarily K\"ahler-Einstein. This partially answers a question of X.X. Chen.
\end{abstract}
\maketitle

\section{Introduction}

One of the major problems in K\"ahler geometry is the study of extremal metrics, especially K\"ahler-Einstein metrics
or constant scalar curvature K\"ahler metrics (cscK). Let $(M,\omega)$ be a compact K\"ahler manifold of complex dimension $n$, and
assume the first Chern class of $M$ to be negative or zero. Then by the celebrated work of Yau \cite{yau1} (see also Aubin \cite{aubin})
we know that $M$ admits a K\"ahler-Einstein metric. The remaining case $c_1(M)>0$ is more difficult and still largely
open. There is a conjecture of Yau \cite{yau2} that the existence of a K\"ahler-Einstein metric should be equivalent
to the stability of the manifold in some algebro-geometric sense.
The behaviour of the Mabuchi energy \cite{mabuchi} plays a central role in this study \cite{tian}.
X.X. Chen and G. Tian in \cite{chentian1} have introduced a family of functionals $E_k$ ($0\leq k\leq n$) which generalize
the Mabuchi energy $\nu_\omega=E_0$, and used $E_1$ in their study of the K\"ahler-Ricci flow. They can also be defined in terms of
the Deligne pairing \cite{phong}, and they should possibly be related to stability.
Bando and Mabuchi showed \cite{bandomabuchi} that K\"ahler-Einstein metrics realize the minimum of the Mabuchi energy,
and recently Song and Weinkove \cite{benjian} showed that K\"ahler-Einstein metrics realize the minimum of the $E_k$,
on the space of metrics with nonnegative Ricci curvature in the anticanonical class for $k\geq 2$, and on the whole 
anticanonical class for $k=1$. Pali \cite{pali} also investigated the question of the lower boundedness of $E_1$.

For the moment we drop the assumption $c_1(M)>0$. Let $P(M,\omega)$ be the set of K\"ahler potentials
for $\omega$, i.e. the set of all smooth real functions $\phi$ such that 
$$\omega_\phi:=\omega+\sqrt{-1}\de\db\phi$$
is positive definite. For $\phi\in P(M,\omega)$ let $\phi_t$ be a smooth path in $P(M,\omega)$ with $\phi_0=0$, $\phi_1=\phi$, and define
\begin{equation*}
\begin{split}
E_{k,\omega} (\phi) &= \frac{k+1}{V} \int_0^1 \int_{M} (\triangle_{\phi_t} \dot\phi_t ) \, \Ric (\omega_{\phi_t})^k \wedge \omega_{\phi_t}^{n-k}dt \\
&- \frac{n-k}{V} \int_0^1 \int_M \dot\phi_t \, ( \Ric (\omega_{\phi_t})^{k+1} - \mu_k \, \omega_{\phi_t}^{k+1}) \wedge \omega_{\phi_t}^{n-k-1}dt,
\end{split}
\end{equation*}
where $\dot\phi_t=\frac{\de}{\de t}\phi_t$, $V=\int_M\omega^n$ is the volume of $M$, and
$$\mu_k=\frac{1}{V}\int_M\Ric(\omega)^{k+1}\wedge\omega^{n-k-1}=(2\pi)^{k+1}\frac{[c_1(M)]^{k+1}\smallsmile[\omega]^{n-k-1}}{[\omega]^n}$$
is a constant that depends only on the cohomology classes $[c_1(M)]$ and $[\omega]$. It can be checked that these functionals do not
depend on the choice of the path. The functional $E_{0,\omega}$ is equal to the Mabuchi energy $\nu_\omega$.
With a simple integration by parts we can rewrite these functionals as
\begin{multline*}
E_{k,\omega}(\phi)\\
=\frac{1}{V} \int_0^1 \int_M \dot\phi_t \left[ (k+1) \binom{n}{k}^{-1} (\triangle_{\phi_t} (\sigma_k(\omega_{\phi_t})) - \sigma_{k+1}(\omega_{\phi_t})) + (n-k)\mu_k \right] \omega_{\phi_t}^n dt,
\end{multline*}
where $\sigma_k (\omega)$ is given by
$$(\omega+t\Ric(\omega))^n = \left( \sum_{k=0}^n \sigma_k(\omega)t^k \right) \omega^n.$$
If at a point $p$ we pick normal coordinates in which the Ricci tensor is diagonal with entries $\lambda_1, \dots, \lambda_n$
then $\sigma_{k}(\omega)$ at $p$ is the $k$-th elementary symmetric function in the $\lambda_i$,
$$\sigma_{k}(\omega) = \sum_{1\leq i_1 < \cdots < i_k\leq n} \lambda_{i_1} \lambda_{i_2} \cdots \lambda_{i_k}.$$ 
For example we have $\sigma_0(\omega)=1$, $\sigma_1(\omega)=R$ the scalar curvature of $\omega$, $\sigma_2(\omega)=\frac{1}{2}(R^2-|\Ric|^2)$.
From this we see that the critical points of $E_{k,\omega}$ are the functions $\phi$ such that $\omega_\phi$ satisfies
\begin{equation}\label{crit}
\sigma_{k+1}(\omega_\phi)-\triangle_{\phi}\sigma_k(\omega_\phi)=\binom{n}{k+1}\mu_k.
\end{equation}
For $k=0$ we see that the critical metrics are precisely cscK metrics. The critical equation for $k=1$ is
\begin{equation}\label{crite1}
R^2-|\Ric|^2-2\triangle R=n(n-1)\mu_1.
\end{equation}

From now on assume that $c_1(M)>0$, that $[\omega]=2\pi[c_1(M)]$ so that $\mu_k=1$ for all $k$. Under these hypotheses Chen and Tian \cite{chentian1} have shown that critical metrics for $E_n$ with positive Ricci curvature must be
necessarily K\"ahler-Einstein. In a recent paper \cite{chen1} X.X. Chen asked the following:

\begin{quest}[X.X. Chen]
Is a critical metric for $E_1$ in the anticanonical class necessarily K\"ahler-Einstein?
\end{quest}

Our result is a partial answer to this question, for all $k$, generalizing Chen and Tian:

\begin{thm}\label{main}
Fix $0\leq k\leq n$ and assume $\omega$ is a K\"ahler metric in the anticanonical class which is critical for
$E_k$ and has nonnegative Ricci curvature. Then $\omega$ is K\"ahler-Einstein. If $k=1$ then the same holds
if $\omega$ satisfies only $R>-n$.
\end{thm}

\begin{rem} G. Maschler \cite{masch} has proved that a critical point of $E_n$ in the anticanonical class
is always K\"ahler-Einstein.
\end{rem}

\section{The proof}

In case $k=0$ the critical metrics are cscK and it is well known that in the anticanonical class
they are necessarily K\"ahler-Einstein, even if the Ricci curvature is not positive (see
for example \cite[Prop. 2.12]{tianbook}). Next we assume that $0<k<n$, that $\Ric(\omega)>0$ and we use the minimum principle.
Let $p\in M$ be a point where $\sigma_k(\omega)$ achieves its minimum,
and pick normal coordinates around $p$ in which the Ricci tensor is diagonal with \emph{positive} entries $\lambda_1,\dots,\lambda_n$.
Notice that all the $\sigma_k(\omega)$ are positive.
Then using the critical equation \eqref{crit} we have that
$$\sigma_{k+1}(\omega)(p)\geq\binom{n}{k+1}.$$
Define the normalized symmetric functions as
$$\Sigma_k(\omega):=\frac{\sigma_k(\omega)}{\binom{n}{k}},$$
so that at $p$ we have $\Sigma_{k+1}(\omega)\geq 1$. Since the $\lambda_i$'s are positive,
a classical inequality (see \cite[Sect. 2.22, Thm. 51]{polya})
says that for every $0\leq k\leq n-1$ we have
\begin{equation}\label{polya}
\Sigma_k^{k+1}(\omega)\geq\Sigma_{k+1}^k(\omega),
\end{equation}
and this holds everywhere on $M$.
Hence we have
$\Sigma_k(\omega)(p)\geq 1$
and since $p$ is the minimum point of $\Sigma_k(\omega)$, we have $\Sigma_k(\omega)\geq 1$ everywhere. Also using \eqref{polya} repeatedly
we have that $$\Sigma_1(\omega)\geq\Sigma_2^{\frac{1}{2}}(\omega)\geq\dots\geq\Sigma_k^{\frac{1}{k}}(\omega)\geq 1$$
everywhere on $M$. Now by hypothesis $\omega\in 2\pi c_1(M)$ so there exists a smooth real function $f$ such that
$\Ric(\omega)=\omega+\sqrt{-1}\de\db f$, and taking the trace of this with respect to $\omega$ we get
$R=n+\triangle f$, which can be rewritten as $\Sigma_1(\omega)=1+\frac{\triangle f}{n}$. Hence $\triangle f\geq 0$ from
which it follows easily that $f$ is constant, and so $\Ric(\omega)=\omega$ and $\omega$ is K\"ahler-Einstein.

If now we only assume that $\Ric(\omega)\geq 0$, we have that $\lambda_i\geq 0$. Using the minimum principle as before,
we see that $\Sigma_{k+1}(\omega)(p)\geq 1$, so it is not zero. From this it follows that $\Sigma_j(\omega)(p)>0$ for $0\leq j\leq k$,
and it is very easy to check that in this case the inequality \eqref{polya} still holds and the reasoning proceeds as before.

If $k=n$ then the critical equation becomes $\triangle\sigma_n(\omega)=0$. Hence $\sigma_n(\omega)=1=\Sigma_n(\omega)$, and
as before $\Sigma_1(\omega)\geq\Sigma_n^\frac{1}{n}(\omega)=1$ and so $\omega$ is K\"ahler-Einstein. G. Maschler called
such metrics \emph{central K\"ahler metrics} and he proved \cite{masch} that the hypothesis $\Ric\geq 0$ can be removed,
using Demailly's holomorphic Morse inequalities.

If $k=1$ then the critical equation is \eqref{crite1} with $\mu_1=1$. In this case the inequality \eqref{polya} holds
for arbitrary $\lambda_i$, and takes the form $|\Ric|^2\geq\frac{R^2}{n}$ (which follows from
$|\Ric-\frac{R}{n}g|^2\geq 0$). With the same argument as above
we get at the minimum $p$ of $R$, $R^2(p)\geq n^2.$
Since we are assuming $R>-n$, we must have $R(p)\geq n$, so $R\geq n$ everywhere and $\omega$ is K\"ahler-Einstein.

\section{Further remarks}

\begin{quest} Is Theorem \ref{main} still true if we do not assume the Ricci curvature to be nonnegative (at least for $k=1$)?
\end{quest}
Notice that this holds for $k=0,n$.

\begin{rem} Further evidence for the validity of Chen's question is provided by this recent result
of Song and Weinkove \cite{benjian}: assuming $M$ has no nonzero holomorphic vector fields, 
the functional $E_1$ is proper in the sense of Tian \cite{tian} if and
only if there exists a K\"ahler-Einstein metric. This leads one to expect that K\"ahler-Einstein metrics
and critical points of $E_1$ coincide.
\end{rem}

\begin{ack}The author thanks Jian Song and Ben Weinkove for many useful discussions, and his advisor
Professor S.-T. Yau for his constant support.
\end{ack}

\end{document}